\documentclass[12pt]{amsart}
\usepackage{amssymb,amscd}
\usepackage[mathscr]{eucal}
\usepackage{pb-diagram}

\newtheorem{thm}{Theorem}

\newtheorem{lem}[thm]{Lemma}

\newtheorem{rem}[thm]{Remark}

\begin{document}

%
%
%

\renewcommand{\o}{\operatorname}
\newcommand{\Spec}{\o{\bold {Spec}}}
\newcommand{\A}{{\mathbb A}}
\newcommand{\N}{{\mathbb N}}
\newcommand{\C}{{\mathbb C}}
\newcommand{\E}{{\o E}}
\newcommand{\Mgn}{{M_{g,n}}}
\newcommand{\CC}{{\mathcal C}}
\renewcommand{\P}{{\mathbb P}}
\newcommand{\bH}{{\mathbb H}}

\newcommand{\M}{{\mathcal M}}
\newcommand{\U}{{\mathcal U}}
\newcommand{\X}{{\mathcal X}}
\newcommand{\bM}{{\bar M}}
\newcommand{\B}{{\mathbb B}}
\newcommand{\MC}{{\mathcal M}_{\mathcal C}}
\renewcommand{\H}{\o{H}}
\newcommand{\RR}{\mathcal{R}}
\newcommand{\R}{{\mathbb R}}
\newcommand{\Z}{{\mathbb Z}}
\newcommand{\Id}{{\mathbb I}}
\newcommand{\Os}{{\mathcal O}}
\newcommand{\zz}{{\mathcal Z}}
\newcommand{\z}{{\zz}^1(M)}
\newcommand{\la}{\langle}
\newcommand{\ra}{\rangle}
\newcommand{\h}{{\mathcal H}^3_{\R}}
\newcommand{\hH}{{\mathscr H}}
\newcommand{\hi}{h^{-1}}
\renewcommand{\Im}{\o{Im}}
\renewcommand{\Re}{\o{Re}}
\newcommand{\Hol}{\o{Hol}(E)}
\newcommand{\Holu}{\o{Hol}_u(E)}
\newcommand{\hol}{\o{hol}}
\newcommand{\her}{\o{Her}}
\newcommand{\Her}{\her(E)}
\newcommand{\fc}{{\mathcal F}(E)} 
\newcommand{\fcu}{{\mathcal F}_u(E)} 
\newcommand{\fcuu}{{\mathcal F}_u(E)} 
\newcommand{\fcl}{{\mathcal F}_l(E)} 
\newcommand{\g}{{\mathfrak g}}
\renewcommand{\k}{{\mathfrak k}}
\newcommand{\gl}{{\mathfrak{gl}}}
\newcommand{\su}{{\mathfrak{su}}}
\newcommand{\ff}{{\mathfrak F}} 
\newcommand{\ga}{{\mathcal G}(E)} 
\newcommand{\gau}{{\mathcal G}_u(E)} 
\newcommand{\gal}{{\mathcal G}_l(E)} 
\newcommand{\Hig}{\o{Higgs}(E)}
\newcommand{\Hom}{\o{Hom}}
\newcommand{\hpg}{\Hom(\pi_1,G)/G}
\newcommand{\hppg}{\Hom(\pi_1,G)\slashslash G}
\newcommand{\hmg}{\Hom(\pi_1,G)}
\newcommand{\Un}{\o{U}(n)}
\newcommand{\Uo}{\o{U}(1)}
\newcommand{\GLn}{\o{GL}(n,\C)}
\newcommand{\GL}{\o{GL}}
\newcommand{\PUoo}{\o{PU}(1,1)}
\newcommand{\PSLtR}{\o{PSL}(2,\R)}
\newcommand{\PSLnR}{\o{PSL}(n,\R)}
\newcommand{\SLnC}{\o{SL}(n,\C)}
\newcommand{\PSLnC}{\o{PSL}(n,\C)}
\newcommand{\SUtwo}{\o{SU}(2)}
\newcommand{\SU}{\o{SU}(2)}
\newcommand{\SLC}{\o{SL}(2,\C)}
\newcommand{\SLR}{\o{SL}(2,\R)}
\newcommand{\Spin}{\o{Spin}(2)}
\newcommand{\Spinm}{\o{Spin}_-(2)}
\newcommand{\Pin}{\o{Pin}(2)}
\renewcommand{\O}{{\mathcal O}}
\renewcommand{\P}{{\Bbb P}}
\newcommand{\Q}{{\Bbb Q}}
\newcommand{\HomU}{\o{Hom}_{\B}(\pi_1(M), \SU)}
\newcommand{\HomR}{\o{Hom}_\B^+(\pi_1(M), \SLR)}
\newcommand{\HomC}{\o{Hom}_\B^+(\pi_1(M), \SLC)}
\newcommand{\HomG}{\o{Hom}_\B^+(\pi_1(M), G)}
\newcommand{\rk}{\o{rank}}
\newcommand{\Map}{\o{Map}}
\newcommand{\db}{\overline{\partial}}
\renewcommand{\d}{\partial}
\newcommand\ka{K\"ahler~}
\newcommand\tM{\tilde{M}}
\newcommand\tS{\tilde{S}}
\newcommand\tX{\tilde{X}}
\newcommand\td{\tilde{d}}
\newcommand\tg{\tilde{g}}
\newcommand\tT{\tilde{T}}
\renewcommand\th{\tilde{h}}
\newcommand\ca{{\mathcal A}}
\newcommand\cb{{\mathcal B}}
\newcommand\cc{{\mathcal C}}
\newcommand\Iso{{\o{Iso}}}
\newcommand\Aut{{\o{Aut}}}
\newcommand\Ad{{\o{Ad}}}
\newcommand\Obj{{\o{Obj}}}
\newcommand\hk{{hyperk\"ahler~}}
\newcommand\slashslash{{/\hspace{-3pt}/}}
\newcommand\Er{E_{\rho}}
\newcommand{\tr}{\o{tr}}

\title[Connections with SU(2)-monodromy]{Explicit connections with SU(2)-monodromy}
\author{Eugene Z. Xia}

\address{
Department of Mathematics\\
National Cheng kung University\\
Tainan 701, Taiwan } \email{ezxia@ncku.edu.tw}
\date{\today}

\begin{abstract}
The pure braid group $\Gamma$ of a quadruply-punctured Riemann
sphere acts on the $\SLC$-moduli $\M$ of the representation variety
of such sphere.  The points in $\M$ are classified into
$\Gamma$-orbits.  We show that, in this case,  the monodromy groups
of many explicit solutions to the Riemann-Hilbert problem are
subgroups of $\SU$.  Most of these solutions are examples of
representations that have dense images in $\SU$, but with finite
$\Gamma$-orbits in $\M$. These examples relate to explicit
immersions of constant mean curvature surfaces.
\end{abstract}
\maketitle


\section{Introduction}
Let $G$ be an algebraic Lie group over $\C$
with Lie algebra $\g$.  Denote by $\O_G$ and $\O_\g$ the categorical
quotients of the $G$-adjoint actions on $G$ and $\g$, respectively.
Let $n \in \N$.
The exponential map $\E : \g \to G$ induces the exponential
maps $\E : \O_\g \to \O_G$ and $\O_\g^n \to \O_G^n$.
Let
$$C = \{ c \in \P^1(\C)^n : c_n = \infty; c_i \neq c_j \text{ if } i
\neq j\},
$$
$$
\X(\g) = \{ X \in \g^n : \sum_{i=1}^n X_i = 0\}.
$$
For each $c \in C$, let $\Sigma = \Sigma(c) = \P^1(\C) \setminus
\{c_1,\cdots,c_n\}$. Fix a base point $p$ and let $\pi_1 =
\pi_1(\Sigma,p)$ be the fundamental group.  The representation
variety $\Hom(\pi_1,G)$ identifies with
$$
\RR(G) = \{A \in G^n : \prod_{i=1}^n A_i = e\}
$$
which has a natural variety structure inherited from $G$.  The
diagonal $G$-adjoint actions on $\X(\g)$ and $\RR(G)$ give
categorical quotients $\U(G)$ and $\M(G)$, respectively. Denote by
$P$ the projections $\RR(G) \to \O_G^n$ and $\X(\g) \to \O_\g^n$.
We do not distinguish $X \in \X(\g)$ and $A \in \RR(G)$ from their
respective images in $\U(G)$ and $\M(G)$, the group is always
assumed to be $G$ and we shorten $\RR(G)$ to $\RR$ and etc, unless
otherwise specified.

The pure braid group $\Gamma$ of $\Sigma$ acts on $\pi_1$, hence, on
$\RR$ and $\M$.  The $\Gamma$-action preserves the fibre of $P$. For
each $a \in \O_G^n$ (resp. $\theta \in \O_\g^n$), let $\RR_a =
P^{-1}(a) \subset \RR$ (resp. $\X_\theta = P^{-1}(\theta)$) and
$\M_a = P^{-1}(a) \subset \M$ (resp. $\U_\theta = P^{-1}(\theta)$).

For $z \in \Sigma(c)$ and $X \in \X$,
$$
D_X = \partial + \sum_{i = 1}^{n-1} \frac{X_i}{z - c_i}dz
$$
is a flat connection on $\Sigma$.   Each $D_X$ induces a
representation $\pi_1 \to G$.  This gives rise to monodromy maps
$\hol : \X \to \RR$ and $\U \to \M. $

The Riemann-Hilbert problem concerns the surjectivity of $\hol$.
It is an existential question, but one may ask the constructive
Riemann-Hilbert question: Given $A \in \RR$, construct $X \in \X$
such that $\hol(X) = A$.
When $n=3$ and $G = \SLC$, there is {\em rigidity}, i.e. $\M_a$
consists of a point for a generic $a \in \O_G^3$. Notice that $\E
\circ \hol$ equals $\E$ on $\U$ and the latter is simple to compute.
If $\M_a$ is not empty, then, up to equivalence of representations,
any $X \in \X$ with $\E(X) = a$ is a solution.

In general, $\M_a$ is a moduli space with positive dimension for a generic
$a$. The constructive Riemann-Hilbert problem involves solutions of
non-linear differential equations and has been solved for the finite
subgroups of $\SU$ with $n = 4$ only recently \cite{Bo1, Bo2, Hi}.
Here we pose a related problem:

\noindent {\bf Question:} Suppose $K < G$. Given $X \in \X$,
determine whether, up to a $G$-inner isomorphism, $\hol(X)(\pi_1) <
K$ or, equivalently, $\hol(X) \in \RR(K)$.

This is an interesting problem in its own right, but also has
applications.  When $G =\SLC$ and $K = \SU$, examples of such $X \in
\X$ are related to explicit immersions of constant mean curvature
surfaces into the Euclidean space ${\mathbb E}^3$, the hyperbolic space  $\bH^3$
and the standard sphere $S^3$ \cite{DPW, RS, Sc, SKKR}.

In a more general context, one may consider $\Sigma$ to be a
punctured Riemann surface of genus $g$.  The mapping class group
$\Gamma$ fixing $\partial \Sigma$ acts on an analogous $\M_a(\SLC)$,
preserving $\M_a(\SU)$, the real points corresponding to $\SU$-representations. Suppose further that $g > 0$. If $\rho \in
\RR_a(\SU)$ and $\rho(\pi_1)$ is dense in $\SU$, then the
$\Gamma$-orbit is dense in $\M_a(\SU)$ \cite{PX1, PX2}.
However this is no longer true when $g = 0$ (our present case)
\cite{PX3}. The results here also provide more such examples in the
case of $g = 0$. \newline

\noindent {\bf Acknowledgement:} The author benefited from
discussions with Philip Boalch while going through \cite{Bo1,Bo2}.

\section{Representation Varieties}
From now on, assume
$n=4$ and $\SU = K < G=\SLC$.  For $H < G$, $A \in \RR$ (resp. $X \in \X$) is said to
be an $H$-class if, up to a $G$-inner isomorphism, $A(\pi_1) < H$ (resp. $\hol(X)(\pi_1) < H$).
For a generic $a$, $\M_a$ is a smooth two dimensional moduli, hence,
no longer rigid. However the points of $\M_a$ are classified into
$\Gamma$-orbits. The finite orbits are actually subvarieties of
$\M_a$ that are rigid in certain sense.
As examples, if $A(\pi_1)$ is finite, then the
$\Gamma$-orbit of $A$ is finite.  The classification of such $A$ is
elaborate, but carried out in \cite{Bo1, Bo2}.

Up to a M\"obius transformation, assume $c = \{0,1,\infty, t\}$. The
space $\O_G^4$ can be identified with $\C^4$ and the projection with
the trace map \cite{BG,Go2}
$$P(A_1,A_2,A_3,A_4) = (\tr(A_1), \tr(A_2), \tr(A_3), \tr(A_4)).$$
Let
$$a = (a_1, a_2, a_2, a_4) \in \O_G^4,$$
$$v = (v_1, v_2, v_3) = (\tr(A_1 A_2), \tr(A_2 A_3), \tr(A_1 A_3)).$$
In this convention, the exponential map $\E : \O_\g \to \O_G$ is
defined by $\E(\theta) = 2 \cos(\pi \theta)$.  By \cite{BG,Go2},
$$\M = \{(a,v) \in \C^7 : f(a,v) = 0\}, \ \  \M_a = \{v \in \C^3 :
f(a,v) = 0\},$$ where
\begin{equation}\label{ideal}
f(a,v) = v_1^2 + v_2^2 + v_3^2 + v_1v_2v_3 - (a_1 a_2 + a_3 a_4)v_1
- (a_1 a_4 + a_2 a_3) v_2
\end{equation}
$$
- (a_1 a_3 + a_2 a_4) v_3 + (a_1^2 + a_2^2 + a_2^2 + a_4^2 + a_1 a_2
a_3 a_4 - 4).
$$

The group $\Gamma$ has three generators $\tau_1, \tau_2$ and
$\tau_3$ with its actions on $\M_a$ \cite{BG, Go2}:

$$
\left[
\begin{array}{c}
v_1 \\
v_2\\
v_3\\
\end{array}
\right] \stackrel{\tau_1}{\longmapsto} \left[
\begin{array}{c}
v_1\\
a_1a_4+a_2a_3-v_1(a_1a_3+a_2a_4-v_1 v_2 - v_3)-v_2\\
a_1a_3+a_2a_4- v_1v_2-v_3\\
\end{array}
\right],
$$
$$\left[
\begin{array}{c}
v_1\\
v_2\\
v_3\\
\end{array}
\right] \stackrel{\tau_2}{\longmapsto} \left[
\begin{array}{c}
a_1a_2+a_3a_4 - v_2v_3-v_1\\
v_2\\
a_2a_4+a_1a_3-v_2(a_1a_2+a_3a_4-v_2v_3-v_1)-v_3\\
\end{array}
\right],
$$
$$\left[
\begin{array}{c}
v_1\\
v_2\\
v_3\\
\end{array}
\right] \stackrel{\tau_3}{\longmapsto} \left[
\begin{array}{c}
a_1a_2 + a_3a_4 - v_3(a_2 a_3+a_1a_4-v_3v_1-v_2)-v_1\\
a_2a_3+a_1a_4- v_3v_1 - v_2\\
v_3\\
\end{array}
\right].
$$

A representation $A \in \RR$ is an $\SU$- or $\SLR$-class if and
only if $a$ and $v$ are real \cite{BG}. Define the interval
$$
I_{s,t} = [\frac{s t - \sqrt{(s^2 - 4) (t^2 - 4)}}{2}, \ \frac{s t +
\sqrt{(s^2 - 4) (t^2 - 4)}}{2}].
$$
Then by \cite{BG},
\begin{lem}\label{su2 connections} A representation $A$ is an
$\SU$-class if and only if $v \in \R^3 \subset \C^3$, $A
\in \M_a$ with $a \in [-2,2]^4 \subset \C^4$ and $I_{a_1,a_2} \cap
I_{a_3, a_4} \neq \emptyset.$
\end{lem}

\begin{rem}\label{rem:R}
Let $a \in \O_G^4$.  Suppose $\Delta < \Gamma$ has finite index. Let
$\M_a^\Delta \subset \M_a$ be the $\Delta$-fixed subvariety. The
subspace $\M_a^\Delta$ is discrete in $\M_a$ and one may determine
whether $\M_a^\Delta$ consists of $\SU$-classes by Lemma \ref{su2
connections} (See \cite{PX3} for an example).
\end{rem}
\section{Flat connections}
Let $\theta \in \O_\g$ and $a = \E(\theta)$ which equals to $2
\cos(\pi \theta)$ in our convention.  Given $A \in \RR_a$, explicit
solutions of $X = (X_1, X_2, X_3, X_4) \in \X_\theta$ with $\hol(X)
= A$ are related to solutions of the Painlev\'e VI equation
\cite{Bo1, Bo2}:

\begin{equation} \label{VI}
\frac{d^2 y}{dt^2 } - \frac{1}{2} (\frac{1}{y} + \frac{1}{y - 1} +
\frac{1}{y - t}) (\frac{d y}{d t})^2 + (\frac{1}{t} + \frac{1}{t -
1} + \frac{1}{y - t})\frac{dy}{dt}
\end{equation}
$$
= \frac{t(y-1)(y-t)}{t^2(t-1)^2}(r_1 + r_2 \frac{t}{y^2} + r_3
\frac{t-1}{(y-1)^2} + r_4 \frac{t(t-1)}{(y-t)^2}),
$$
where
$$
r_1 = \frac{(\theta_4 - 1)^2}{2}, \ r_2 = -\frac{\theta_1^2}{2}, \
r_3 = \frac{\theta_3^2}{2},\ r_4 = \frac{(1-\theta_2^2)}{2}.
$$

The group $\Gamma$ acts on the moduli of flat connections, hence, on
the solutions of the Painlev\'e VI equation \cite{Bo1, Bo2}.  An
algebraic solution to the Painlev\'e VI equation has a finite
$\Gamma$-orbit (See \cite{Bo1, Bo2}) with isotropy subgroup $\Delta
< \Gamma$. Once such a solution is found, one may deform it in many
ways to obtain $\Delta$-fixed families of solutions. For
example, equation (\ref{VI}) is actually a family of equations
parameterized by $\theta$. Hence if $y$ is an explicit solution,
then $y$ is a solution to a family of Painlev\'e VI equations if
$\theta$ is deformed in such a way that the right hand side of
Equation (\ref{VI}) remains constant. This family of equations then
correspond to a $\theta$-family of solutions $\Theta \subset
\X^\Delta$.  Set $\Theta_\theta  = \Theta \cap \X_\theta$. It then
follows that

\begin{thm}\label{main}
$\hol(\Theta_\theta) \subset \M_a^\Delta$.  Hence if $\M_a^\Delta$
consists of only $\SLR$- or $\SU$-classes (See Remark \ref{rem:R}),
then $\Theta_\theta$ consists of only $\SLR$- or $\SU$-classes,
respectively.
\end{thm}

\section{Examples}
Consider (\cite{Bo2} \S 3,
Example 3). Let $A = \hol(X)$. Then $A(\pi_1) < K$ is contained in
the symmetry group of the tetrahedron and the $\Gamma$-orbit of $A$
consists of exactly two points. By a direct computation, $\Theta$ is
parameterized by the affine variety
$$
\{\theta \in \C^4 : -\theta_2^2 + \theta_3^2 = 0, 1 - \theta_1^2 - 2
\theta_4 + \theta_4^2 = 0\}.
$$
Moreover $A  = (0,1,0) \in \M_{(1, -1, -1, -1)}$ and
\begin{equation}\label{two point orbit}
\Delta = \langle \tau_2, \tau_1^2, \tau_3^2 \rangle.
\end{equation}
By a direct computation, using Gr\"obner bases, the
subvariety $\M^\Delta \subset \M$ is defined by the ideal
$$
(-4 v_3 + v_2^2 v_3, -2 v_1 - v_2 v_3, 4 - 2 a_3^2 - 2 a_4^2 + a_3^2
a_4^2 + a_3^2 v_2 - a_4^2 v_2 - v_2^2, a_2 - a_3, a_1 + a_4)
$$
Hence if $a = (a_1, a_2, a_2, -a_1)$, then
$$\M_a^\Delta = \{(0, 2 - a_1^2, 0), (0, a_2^2 - 2, 0)\}.$$
Furthermore, if $a \in [-2,2]^4$ satisfies the additional hypothesis
of Lemma~\ref{su2 connections}, then $\M_a^\Delta$ consists of
$\SU$-classes.
\begin{rem}
The $\SU$-classes in Theorem \ref{main} are the ones found in
\cite{PX3}.  We emphasize here that the matrices in $\Theta$ can be
explicitly computed, but since the formulas are rather complicated,
we refer to (\cite{Bo2}, Appendix A) for details.
\end{rem}
For the case of (\cite{Bo2} \S 3, Example 4), $a = (-2 + a_3^2, a_3,
a_3, -1)$ and
$$
\M_a^\Delta = \{(0,1,0), (0,1, -3a_3 + a_3^3), (-3a_3 + a_3^3, 1, 0)
\}.
$$
For (\cite{Bo2} \S 3, Example 5), $a = (a_3, a_3, a_3, 0)$ and
$$
\M_a^\Delta = \{(1,-2 + a_3^2,1), (1,1,1), (-2 + a_3^2, 1, 1), (1,
1, -2 + a_3^2) \}. $$ Case (\cite{Bo2} \S 3, Example 6) is rigid.

%
%
%
%
%
%
%
%
%
%
%
%

One can similarly work out the octahedron cases \cite{Bo2} and the
icosahedron cases \cite{Bo1}.
\section{Conclusions}
To summarize the algorithm:  For
$[\Gamma : \Delta] < \infty$, let $\M^\Delta \subset \M$ be the
$\Delta$-fixed subvariety. Apply the methods in \cite{Bo1, Bo2,
Hi} to compute the family $\Theta \subset \X^\Delta$.  Let $\theta
\in P(\Theta)$ and $a = \E(\theta) \in \O_G^4$. The subspace
$\M_a^\Delta$ is discrete in $\M_a$ and one may determine whether
$\M_a^\Delta$ consists of $\SLR$- or $\SU$-classes by Lemma \ref{su2
connections}.
When $n = 3$, $\M_a$ consists of a single point for a generic $a$, hence, trivially fixed by $\Gamma$. In this sense, the points in $\M_a^\Delta$ may be thought of as rigid.




\end{document}